\documentclass[]{interact}

\usepackage{hyperref}

\usepackage{latexsym} 
\usepackage{graphics}
\usepackage{epsfig}
\usepackage{graphicx}
\usepackage{tikz-qtree}
\usepackage{tikz}
\usepackage{amsmath}
\usepackage{amssymb}
\usepackage{dirtytalk} 
\usepackage{bbm} 

\usepackage{pgfplots}
\usepgfplotslibrary{groupplots}
\usepackage[utf8]{inputenc}
\usepackage[english]{babel}

\pgfplotsset{width=11cm,compat=newest}
\usepackage{algorithm}
\usepackage[noend]{algpseudocode}
\makeatletter
\def\BState{\State\hskip-\ALG@thistlm}
\makeatother
\algdef{SE}[DOWHILE]{Do}{doWhile}{\algorithmicdo}[1]{\algorithmicwhile\ #1}%

\usepackage{varwidth}
\usepackage{booktabs} 
\usepackage{epsfig}
\usepackage{amsfonts}
\usepackage{array}

\usepackage{rotating}
\usepackage{tabularx}
\usepackage{color}
\usepackage{multirow}
\usepackage{url}
\usepackage{lscape}
\usepackage{xargs} 
\newcolumntype{x}[1]{>{\centering\let\newline\\\arraybackslash\hspace{0pt}}p{#1}}

\usepackage{epstopdf}
\usepackage{subfigure}

\usepackage{natbib}
\bibpunct[, ]{(}{)}{;}{a}{}{,}
\renewcommand\bibfont{\fontsize{10}{12}\selectfont}

\theoremstyle{plain}
\newtheorem{theorem}{Theorem}[section]
\newtheorem{lemma}[theorem]{Lemma}

\theoremstyle{definition}
\newtheorem{definition}[theorem]{Definition}

\theoremstyle{remark}

\begin{document}

\articletype{}
\title{ Stochastic Dynamic Programming Heuristic for the $(R,s,S)$ Policy Parameters Computation}

\author{
\name{A. Visentin\textsuperscript{a}\thanks{CONTACT A.Visentin. Email: andrea.visentin@ucc.ie}, S. Prestwich\textsuperscript{a}, R. Rossi\textsuperscript{b}, S. A. Tarim\textsuperscript{c}}
\affil{\textsuperscript{a}School of Computer Science \& IT, University College Cork, Ireland; \textsuperscript{b}University of Edinburgh Business School, Edinburgh, UK; \textsuperscript{c}Cork University Business School, University College Cork, Ireland}
}

\maketitle

\begin{abstract}
The $(R,s,S)$ is a stochastic inventory control policy widely used by practitioners. In an inventory system managed according to this policy, the inventory is reviewed at instant \textit{R}; if the observed inventory position is lower than the reorder level \textit{s} an order is placed. The order's quantity is set to raise the inventory position to the order-up-to-level \textit{S}.
This paper introduces a new stochastic dynamic program (SDP) based heuristic to compute the $(R,s, S)$ policy parameters for the non-stationary stochastic lot-sizing problem with backlogging of the excessive demand, fixed order and review costs, and linear holding and penalty costs. In a recent work, \cite{visentin2021computing} present an approach to compute optimal policy parameters under these assumptions. 
Our model combines a greedy relaxation of the problem with a modified version of Scarf's $(s,S)$ SDP. A simple implementation of the model requires a prohibitive computational effort to compute the parameters. However, we can speed up the computations by using K-convexity property and memorisation techniques. The resulting algorithm is considerably faster than the state-of-the-art, extending its adoptability by practitioners.
An extensive computational study compares our approach with the algorithms available in the literature.
\end{abstract}

\begin{keywords}
Inventory; demand uncertainty; stochastic lot-sizing; dynamic programming; control policy
\end{keywords}

\section{Introduction}
The computation of solutions for the non-stationary stochastic lot-sizing problem is a well-developed branch of inventory control. The stochasticity of the demand allows modelling the uncertainty of real-world problems, while its non-stationarity allows considering seasonality or life cycle of products. Under this setting, the inventory must satisfy a demand represented by a set of stochastic variables of known probability distributions generally considered independent. \cite{arrow1951optimal} is considered to be the first known work on stochastic inventory models.

The single-item, single-echelon, non-stationary lot-sizing under ordering, holding and penalty cost is an important class of inventory problems. The problem considers a time horizon split in $T$ periods. A wide variety of policies has been developed to manage these systems (\cite{silver1981operations}). A policy defines when an order has to be placed and its quantity. According to the classification proposed in \cite{bookbinder1988strategies}, three policies have been used to deal with stochasticity: static, dynamic and static-dynamic strategies. 

In the static uncertainty, also known as $(R, Q)$, the timing and quantity of the orders are fixed at the beginning of the time horizon; this policy does not react to demand realisations. 

The dynamic strategy, $(s,S)$ policy, checks the inventory level in each period. If the inventory is lower than the order level $s$ an order to raise it to the order-up-to-level $S$ is placed. This policy allows to react to unforeseen demand realisations quickly,  \cite{scarf1959optimality} proves its optimality when the review costs are ignored. However, this policy suffers from a high degree of setup-oriented nervousness (\cite{tunc2011cost,de1997nervousness}); meaning that the order timings frequently change, limiting its practical applicability. 

The static-dynamic strategy, $ (R, S)$ policy, aims to tackle this issue by fixing the replenishment times at the beginning of the time horizon. In this policy, an order that raises the inventory up to $S$ is placed every $R$ periods. Knowing the replenishing time in advance allows to deal better prices and schedule joint deliveries. 

When the demand is non-stationary or the time horizon is finite these policy parameters vary across the planning horizon assuming the $(R_t, Q_t)$, $(s_t, S_t)$ and $(R_t, S_t)$ form, for $t \in [1, \dots, T]$.

The $(R,s,S)$ policy is a generalisation of the dynamic and static-dynamic strategies. In the $(R, s, S)$ policy, the inventory level is assessed at review intervals $R$ if it falls under the $s$ level, an order is placed; the order raises the inventory level to $S$.  If the cost of reviewing the inventory is null, the policy reviews the inventory in each period behaving as the $(s, S)$ one. If the $s$ level is set equal to $S$, an order is placed at each review. In the non-stationary stochastic problem configuration, the policy parameters change across the time horizon, assuming the $(\textbf{R}, \textbf{s}, \textbf{S})$ form.

The $(R, s, S)$ is widely used by practitioners (\cite{silver1981operations}). In the case of stochastic non-stationary problems, three sets of parameters have to be jointly optimised to minimise the expected cost. This task has been considered extremely difficult. In a recent work, \cite{visentin2021computing} introduce the first algorithm to compute the optimal policy parameters. They apply a branch-and-bound approach to explore the possible replenishment plans while computing the order levels and order-up-to-levels using stochastic dynamic programming (SDP). While their method computes the optimal set of parameters, it struggles to scale to big problems, limiting its applicability by practitioners.
In this work, we fill this gap in the literature by:
\begin{itemize}
\item
presenting a relaxation that allows a greedy computation of the replenishment cycles. We combine it with an SDP formulation for the computation of $(R,s,S)$ policy parameters;
\item
introducing computational enhancements that make the model computable in reasonable time;
\item
analysing an extensive numerical study that shows that the heuristic computational effort significantly outperform the optimal method;
\item 
investigating the problem configurations for which the policy computed by the heuristic differs from the optimal one.
\end{itemize}

The paper is structured as follows. A survey of the literature is presented in section
\ref{sec:literature_review}. Section
\ref{sec:prob_description} provides the description of the problem and of the best-known solution, later used as a comparison. Section
\ref{sec:model} introduces the relaxation used, the greedy approach and the computational enhancements. Section
\ref{sec:experimental} shows a comprehensive numerical study.
Finally, Section \ref{sec:conclusions} concludes the paper.

\section{Literature review} \label{sec:literature_review}
This section surveys the relevant stochastic lot-sizing literature. In the first part, we position our approach in comparison to other inventory control policies. We then analyse recent practical applications of the $(R, s, S)$ policy.

An inventory control policy defines when: to assess the inventory, place an order, and the size of the order. The problem of computing policy parameters to satisfy a stochastic demand appears in a wide variety of industrial settings, and it has been extensively investigated in the literature \cite{silver1981operations}. \cite{bookbinder1988strategies} propose a broad framework of inventory control strategies: static uncertainty ($(R, Q)$ policy), dynamic uncertainty ($(s,S)$ policy), and static-dynamic uncertainty ($(R,S)$ policy). It classifies the approaches based on when the replenishments' decisions are taken, if at the beginning of the planning horizon or after realising a period demand. In the $(R,Q)$ policy, the full replenishment plan is fixed at the beginning. A fixed ordering plan is preferred in industrial settings where rigid production/shipment plans are needed. For these reasons, the computation of this policy under uncertainty has been widely investigated, e.g. \cite{sox1997dynamic,meistering2017stabilized,tunc2021mixed}.
\cite{scarf1959optimality} proves that the $(s,S)$ policy (dynamic strategy) is cost-optimal. In this policy, the decision to place an order and its quantity are taken after observing the demand. This policy is particularly effective in dealing with unexpected demand realisations. Recent works involving this policy are \cite{jiao2017stochastic,xiang2018computing,azoury2020optimal}. The static-dynamic uncertainty ($(R,S)$ policy) fixes the replenishment moments at the beginning of the planning horizon and decides the size when placing the order. This policy is preferred because it reduces the setup-oriented nervousness \cite{tunc2011cost}, a known order schedule also allows better deals with the carriers. We refer the readers to relevant studies on this policy, e.g. \cite{tarim2004stochastic,rossi2015piecewise, tunc2018extended}. \cite{ma2019stochastic} presents a survey of stochastic inventory control policy computation.
However, \cite{bookbinder1988strategies} classification does not take into account the stock-taking cost commonly present in real-world problems, e.g. \cite{fathoni2019development,christou2020fast}. The $(R,s,S)$ policy has a lower expected cost compared to the $(s,S)$ one when a cost for assessing the inventory is considered. As mentioned in the introduction, the $(R,s,S)$ can be seen as a generalisation of the $(s, S)$ and $(R,S)$ policies.

The $(R, s, S)$ policy has a vast number of applications in the literature; due to a reduced nervousness compared to the $(s,S)$ and a better cost-performance than the $(R,S)$.
These policies have also been studied for different problem
configurations. \cite{schneider1991empirical,schneider1995power} introduce two heuristics to compute $(R,s,S)$ parameters in a two-echelon inventory system with one warehouse and multiple retailers.
\cite{strijbosch2002simulating} propose a technique to simulate an $(R,s,S)$ inventory system in which the parameters remain constant.
It can compute fill rates or find parameters values to achieve a prescribed service level.  \cite{chen2009coordinated} adopt a hedge-based $(R,s,S)$ policy portfolio, with constant parameters in the short term, for a multi-product inventory control problem. In \cite{cabrera2013solving} the $(R,s,S)$ policy is used to manage the inventory of multiple warehouses. \cite{goccken2015r} use a simulation optimisation technique to determine the optimal policy for distribution centres in a two-echelon inventory system with lost sales. \cite{johansson2020controlling} use an $(R,s,S)$ policy for controlling one-warehouse, multiple-retailer inventory systems; their configuration is motivated by a real problem faced by a company selling metal sheet products. 
In the surveyed papers, the policy parameters are optimised and kept constant across the time horizon, or the $R$ value is a given fixed value reducing the problem to an $(s,S)$ policy. For example, \cite{lagodimos2012computing} solves the continuous-time problem with stationary demand;  \cite{christou2020fast} extend their work to consider the order quantity as a multiple of given batch size. This is due to the complexity in jointly optimising the three sets of parameters. The additional cost of using a stationary policy when the demand varies is well known \cite{tunc2011cost}.\cite{visentin2021computing} introduce the first optimal approach to compute $(R,s,S)$ policy parameters with stochastic non-stationary demand. However, their approach requires considerable effort to solve big instances, limiting its usability for practitioners. 

The survey presented in this section places our work in the stochastic lot-sizing literature. The $(R,s,S)$ policy has a wide variety of applications due to clear advantages over other policies. The algorithm presented herein aims to boost its adoption by providing a heuristic that computes near-optimal policies using a fraction of the computational effort compared to the state-of-the-art.

\section{Problem description} \label{sec:prob_description}
This work considers the single-item, single-stocking location, stochastic
inventory control problem over a $T$-period planning horizon. 
The $(R,s,S)$ policy defines three aspects of inventory management: the timing of inventory reviews, when an order is placed, and the order's size. A review takes place when the inventory level in the warehouse is assessed; these moments are fixed at the beginning of the time horizon. An order can only be placed after a review takes place. The interval between two review moments represents a replenishment cycle.

The demand's stochasticity and non-stationarity of period $t$ are modelled through the random variable $d_t$. Demands are independent variables with a known probability distribution. Cumulative demand of periods $t$ to the beginning of period $j$ takes the form of $d_{t,j}$ with $j>t$.  If the demand in a given period exceeds the on-hand inventory, the excess is backlogged and carried to the next period. In Section \ref{sec:unit_cost}, we extend the model to the lost-sales configuration, where the exceeding demand is lost; a common approach when competitors' products are available. Under these assumptions, the $(R,s,S)$ policy takes the vectorial form form $(\textbf{R},\textbf{s},\textbf{S})$, with $\textbf{R}=(R_1,\dots,R_T)$; where $R_t$ , $s_t$ and $S_t$ denote respectively the length, the reorder-level and order-up-to-level associated with the $t$-th inventory review. 

Policies are compared based on their expected cost. Stocktaking has a fixed cost of $W$. We denote by $Q_t$ the quantity of the order placed in period $t$. Ordering costs are represented by a
fixed value $K$ and a linear cost, but we shall assume that the variable cost is zero without loss of generality. The extension of our solution to the case of a variable production/purchasing cost is straightforward, as this cost can be reduced to a function of the expected closing inventory level at the final period \cite{tarim2004stochastic}. At the end of each period, a holding cost $h$ is charged for every unit carried from one period to the next. In case of a stockout, a penalty cost $b$ is charged for each item and period. We denote with $I_t$ the closing inventory level for period $t$, making $I_0$ the initial inventory.

We consider the problem of computing the $(\textbf{R},\textbf{s},\textbf{S})$ policy parameters that minimize the expected total cost over the planning horizon. The order quantity $Q_t$  is fixed at every review moment before the demand realisation using:
\begin{equation}
Q_t \triangleq q_t(S_t, s_t, I_{t-1})  \triangleq \left\{\begin{matrix}
S_t - I_{t-1} & \text{if }(t\text{ is a review period}) \wedge (I_{t-1} < s_t)\\ 
0 & else
\end{matrix}\right.
\end{equation}
the order is placed only if $t$ is a review period and the open inventory is below the order level $s_t$. For the sake of brevity, in the following formulas we use $Q_t$ as a replacement for $q_t(S_t, s_t, I_{t-1})$.

The problem of computing the optimal $(\textbf{R},\textbf{s},\textbf{S})$ can be formulated as follow:
\begin{equation}
\label{eq:problem_eq}
C_1(I_0) \triangleq  \underset{(\textbf{R},\textbf{s},\textbf{S})}{\min} f_{1}(I_0, Q_1, R_1) + E[ C_{1+R_1}(I_0+ Q_1 - d_{1,1+R_1 }) ]   
\end{equation}

Where $C_1(I_0)$ is the expected cost of the optimal policy parameters starting at period $1$ with the initial inventory $I_0$. In general, $C_t(I_{t-1})$ represent the expected inventory cost of starting at period $t$ with open inventory $I_{t-1}$. While, $f_{t}(I_{t-1}, Q_t, R_t)$ is the expected cost of a review cycle starting in period $t$ and ending up in period $t+R_t$; it comprises review, ordering, holding and penalty cost for the review cycle. 

\begin{align}
\label{eq:immediatecost_rss}
f_{t}(I_{t-1}, Q_t, R_t) &\triangleq & K\mathbbm{1}\{ Q_t>0 \}  + W + \sum_{i = 1}^{R_t}E[h \max(I_{t-1}-d_{t,t+i} + Q_t, 0)\nonumber \\  
&&+ b \max(-I_{t-1}-Q_t+ d_{t,t+i}, 0)]
\end{align}

$C_t(I_{t-1})$ values can be computed recursively when all the policy parameters are computed using the following formula:
\begin{align}
\label{eq:recursion_rss}
C_t(I_{t-1}) \triangleq  f_{t}(I_{t-1}, Q_t, R_t) + E[ C_{t+R_t}(I_{t-1} + Q_t - d_{t,t+R_t }) ]  ) 
\end{align}
until the base case is reached:
\begin{align}
C_{T+1}(I_{T}) \triangleq  0
\end{align}
For a given $(\textbf{R},\textbf{s},\textbf{S})$ parameters set, this formulation allows to compute the expected policy cost. However, the number of combinations of parameters is exponential, making this approach unusable for the computation of optimal ones.

\subsection{Branch-and-bound approach}
\cite{visentin2021computing} present the first algorithm for computing the optimal parameters for the $(\textbf{R},\textbf{s},\textbf{S})$ problem. Their work is based on the following lemma:
\begin{lemma}
If the replenishment cycles ($\textbf{R}$) are fixed, the problem is reduced to a particular version of the $(\textbf{s}, \textbf{S})$ policy computation and can be solved to optimality using Scarf's SDP (\cite{scarf1959optimality}). 
\end{lemma}

In this case, the problem is formulated as:
\begin{equation}
\label{eq:problem_eq_ss}
\widehat{C}_1(I_0) =  \underset{\textbf{R}}{\min} f_{1}(I_0, Q_1, R_1) + E[ C_{1+R_1}(I_0+ Q_1 - d_{1,1+R_1 }) ]  
\end{equation}
Where $\textbf{s}$ and $\textbf{S}$ are dependent on $\textbf{R}$. The proposed baseline compute the optimal replenishment cycles by testing all $\textbf{R}$ possible combinations and computing the optimal policy cost for each of them. Their best technique, our comparison in the experimental section, uses BnB to avoid recomputations and prune sub-optimal $\textbf{R}$ assignment. Optimal $s_t$ and $S_t$ levels can be computed by considering only future periods when $R_t$ is fixed, ignoring the expected opening inventory level; this is not valid for the computation of the $\textbf{R}$ vector.

\section{Heuristic technique} \label{sec:model}
The heuristic introduced in this work aims to compute locally optimal $R_t$ values to produce a near-optimal $(\textbf{R}, \textbf{s}, \textbf{S})$ policy. The main idea is to move the assignment of the decision variable $R_t$ at period $t$ and do not fix all of them at the beginning of the time horizon such as in Equation \ref{eq:problem_eq_ss}. This can be done by transforming the recursive Equation \ref{eq:recursion_rss} into:
\begin{align}
\label{eq:recursion_rss_local}
\widehat{C}_t(I_{t-1}) =  \underset{R_t}{\min} f_{t}(I_{t-1}, Q_t, R_t) + E[ C_{t+R_t}(I_{t-1} + Q_t - d_{t,t+R_t }) ]  ) 
\end{align}
Solving this recursion could lead to different optimal $R_t$ for different opening inventory levels $I_{t-1}$. For example, if the opening inventory level is slightly higher than $s_t$ but considerably lower than $S_t$ an order is not placed, but the next review cycle might be shortened. However, in the $(R,s,S)$ policy, the review cycles are fixed at the beginning of the time horizon and not after the demand realisation. This is the reason why we need to know the probability of the opening inventory level to determine the optimal $R_t$. 

Our heuristics consists of choosing a locally optimal $R_t$ assuming that an order is placed in period $t$ and the possibility of placing a negative order. We define these locally optimal replenishment cycles as $R^a_t$. The independence of the replenishment cycles is similar to the $(R,S)$ policy, and the negative order relaxation is widely used, e.g. \cite{ozen2012static}.

Knowing the expected cost of future periods $\widehat{C}_j$ with $j>t$, it is possible to compute the optimal $s_t$ and $S_t$ for that specific replenishment cycle $R_t$ using SDP. The best $S_t$ is the value that minimizes $\widehat{C}_t(S_t)$, since we place an order to reach the point with the lowest future expected cost.
\begin{align}
\label{eq:scomputation}
S_t =  \arg \min_{I_{t-1}} \widehat{C}_t(I_{t-1})
\end{align}
So, assuming that an order is placed, the best replenishment cycle is the one that has the lowest cost after the inventory level is topped up to $S_t$:
\begin{align}
\label{eq:rheuristic}
R^a_t \triangleq  \arg \min_{R_t}  \widehat{C}_t(S_t)
\end{align}
and 

As mentioned above, the computation of $\widehat{C}_t$ requires the expected costs of future periods $\widehat{C}_j$ with $j>t$, which are dependent on the optimal $R_j$. We relaxed the cost function by defining $C^a_t$ as the expected cost of using local optimal $R^a_j$ for all periods $j$ after $t$. Given $C^a_{T+1}(I_{T}) =  0$, it is possible to compute the relaxed cost function in a backward way using the following approximate SDP functional equation:
\begin{align}
\label{eq:recursion_rss_edit}
C^a_t(I_{t-1}) \triangleq f_{t}(I_{t-1}, Q_t, R^a_t) + E[ C^a_{t+R^a_t}(I_{t-1} + Q_t - d_{t,t+R^a_t }) ]
\end{align}

     This formula computes a near-optimal replenishment schedule $\textbf{R}^a$, and the set of order and order-up-to levels optimal for that given schedule. Due to the relaxation, $\textbf{R}^a$ can differ from the optimal $\textbf{R}$; however, as the experimental section shows, this event is rare.

The resulting approximate SDP formulation is more complex than the $(s,S)$ one, making the computational effort required to solve it prohibitive. This is mainly due to the computation of the expected cycle cost (Equation \ref{eq:immediatecost_rss}); its computation involves three variables in each period: current inventory, order size and length of the replenishment cycle. This computational effort can be considerably reduced applying the K-convexity property (\cite{scarf1959optimality}) used in the $(s,S)$ SDP formulation. The deployment of search reduction and memoisation techniques further improve the performances, and it has a crucial impact on the applicability of this model. In the next subsections, we present the pseudocode for the solution and how these enhancements affect it.

\subsection{Pseudocode}
Algorithm \ref{alg:sdp} shows the procedure to compute the heuristics backwards. Lines 1-2 contains the boundary condition. Line 3 goes through all the periods in a backwards order. Line 5 searches through all the possible replenishment cycles, line 6 through all the inventory levels and line 7 through all the possible order quantities. Lines 12-13 save the current value of $R_t^a$ according to Equation \ref{eq:rheuristic}, while line 14 updates the relative expected costs, Equation \ref{eq:recursion_rss_edit}.
\begin{algorithm}
\caption{RsS-SDP()}\label{alg:sdp}
\begin{algorithmic}[1]
\For{$i$ from $min\_inventory$ to $max\_inventory$}
    \State $C^a_{T+1}(i) = 0$
\EndFor
\For{$t$ from $T$ down to $1$} 
    \State $best\_review\_cost \gets \infty$
    \For{$r$ from $1$ to $T-t +1$}
        \For{$i$ from $min\_inventory$ to $max\_inventory$}
            \State $C_{cycle}(i) \gets \infty$
            \For{$q$ from $0$ to $max\_order$}
                \State $expected\_cost \gets f_t(i,q,r) + E[ C^a_{t+r}(i + q - d_{t,t+r}) ]$
                \If{$expected\_cost < C_{cycle}(i)$}
                    \State $C_{cycle}(i) \gets expected\_cost$
                \EndIf
            \EndFor
        \EndFor
        \If{$min(C_{cycle}) < best\_review\_cost$}
            \State $R^a_t \gets r$
            \State $C_{t}^a \gets C_{cycle}$
            \State  $best\_review\_cost \gets min(C_{cycle})$
        \EndIf
    \EndFor
\EndFor
\end{algorithmic}
\end{algorithm}

For clarity and for the sake of the enhancements, we separate the computation of the immediate cost. Let $\zeta_{t,t+j}$ be a value of the random variable $d_{t,t+j}$ and $P(\zeta_{t,t+j})$ be the probability of assuming that value. Algorithm \ref{alg:sdprscost} computes the immediate cost, Equation \ref{eq:immediatecost_rss}. 
\begin{algorithm}
\caption{$f_t(i,q,r)$}\label{alg:sdprscost}
\begin{algorithmic}[1]
\State $cost \gets W $
\If{$q > 0 $}
    \State $cost \gets cost + K $
\EndIf
\For{$j$ from $1$ to $r$} 
    \For{each $\zeta_{t,t+j}$ value of $d_{t,t+j}$}
        \State $close\_inv \gets i + q - \zeta_{t,t+j}$
        \If{$close\_inv \geq 0$}
            \State $cost \gets cost + h \: close\_inv \:P(\zeta_{t,t+j})$
        \Else
            \State $cost \gets cost - b \: close\_inv \:P(\zeta_{t,t+j})$
        \EndIf
    \EndFor
\EndFor
\Return $cost$
\end{algorithmic}
\end{algorithm}

\subsection{K-convexity}\label{sec:kconv}

We can exploit the property of K-convexity presented in \cite{scarf1959optimality} in solving the dynamic program. This approach is widely used to optimise the $(s,S)$ SDP computation.

The property is defined as:
\begin{definition}
Let $K\geq0$, then function $f(x)$ is K-convex if:
\begin{equation} \nonumber
    K + f(a + x) - f(x) - a \left( \frac{f(x)-f(x-b)}{b} \right) \geq 0
\end{equation}
for all positive $a$, $b$ and $x$.
\end{definition}

\cite{scarf1959optimality} shows that considering $s_t^*$ and $S_t^*$ the optimal reorder level and order up-to level for period $t$:
\begin{equation}\label{eq:kconvexity}
    C_t(I_{t-1} )=  \left \{ \begin{matrix}
f(I_{t-1},0) + E[C_{t+1} (I_{t-1} -d_t)] & s_t^* \leq I_{t-1} \leq S_t^* \\ 
f(I_{t-1},0) + E[C_{t+1} (S_t^* -d_t)] + K & 0 \leq I_{t-1} < s_t^*
\end{matrix} \right.
\end{equation}
This is done by computing the $C_t(y)$ for different values of $y$ starting from an upper bound of $S_t$. The value $y$ is then decremented, and the lowest value of $C_t$ is remembered. When the cost is greater than $C_t + K$ the search terminates. $S_t$ is the inventory level in which the cost assumes the minimum value, $s_t$ is the one in which we stop the search. This approach greatly speeds up the computation of the SDP.

Similarly to the computation of the $(s, S)$ policy, we can use the K-convexity property for the $(R,s,S)$. Considering the Equation  \ref{eq:recursion_rss_edit}, for a fixed $R_t$ the problem is reduced to an $(s,S)$ one with the next $R_t -1$ periods in which an order can not be placed.

\begin{algorithm}
\caption{RsS-SDP-KConv()}\label{alg:sdprsskconv}
\begin{algorithmic}[1]

\For{$i$ from $min\_inventory$ to $max\_inventory$}
    \State $C^a_{T+1}(i) = 0$
\EndFor
\For{$t$ from $T$ down to $1$} 
    \State $best\_review\_cost \gets \infty$
    \For{$r$ from $1$ to $T-t +1$}
        \State $best\_cycle\_cost \gets \infty$
        \For{$i$ from $max\_inventory$ down to $min\_inventory$}
            \State $ C_{cycle}(i) \gets f_t(i,0,r) + E[ C^a_{t+1}(I_{t-1}+ Q_t - d_t) ]$
            \If{$C_{cycle}(i) < best\_cycle\_cost$}
                    \State $best\_cycle\_cost \gets C_{temp}(i)$
                    \State $S_t^{cycle} \gets i$
            \EndIf
            \If{$C_t^{cycle}(i) > best\_cycle\_cost + K$}
                    \State $s_t^{cycle} \gets i$
                    \State \textbf{break for}
            \EndIf
        \EndFor
        \If{$best\_cycle\_cost < best\_review\_cost$}
            \State $R^a_t = r$
            \State $best\_review\_cost \gets best\_cost\_cycle$
            \For{$i$ from $min\_inventory$ to $s_t^{cycle}$}
                \State $C_t^{a}(i) \gets C_t^{cycle}(s_t)$
            \EndFor
            \State $C_{t}^a \gets C_{cycle}$

        \EndIf
    \EndFor
\EndFor
\end{algorithmic}
\end{algorithm}

Algorithm \ref{alg:sdprsskconv} shows the pseudocode of the enhanced SDP, clarifying the improvement's reason. For a fixed review cycle length $R_t$, there is no need to search for the best order quantity $Q_t$. When the order level $s_t$ is determined, the lower inventory levels assume the same expected cost.

\subsection{Cycle Cost Memoisation}\label{sec:memoisation}
The calculation of the cycle cost is particularly time demanding. There is a summation of expected costs over multiple periods. However, it is possible to identify situations in which the same computations occur multiple times. 
 Let $l_t(I_t, R_t)$ be the function that computes the holding and penalty expected cost of starting at the end of period $t$ with closing inventory $I_t$ and with the next review moment in $R_t$ periods. This new function is defined as:
\begin{equation} \label{eq:l}
    l_t(I_t, R_t) \triangleq \sum_{i = 1}^{R_t}E[h \max(I_{t}- d_{t+1,t+i} , 0)+ b \max(-I_{t} + d_{t+1,t+i}, 0)]
\end{equation}
considering $d_{i,j} = 0$ when $i = j$.
Equation \ref{eq:immediatecost_rss} can be rewritten as:
\begin{equation}
     f_{t}(I_{t-1}, Q_t, R_t) = K\mathbbm{1}\{ Q_t>0 \}  + W + l_t(I_{t-1}-d_{t,t+i} + Q_t,R_t)
\end{equation}
The $l_t(I_t, R_t)$ function can be computed in a recursive way:
\begin{equation} \label{eq:recursiveimmediatecost_rs}
    l_t(I_t, R_t) = h \max(I_{t} , 0) + b \max(-I_{t}, 0) + E[l_{t+1}(I_{t} - d_{t+1}, R_t - 1)]
\end{equation}
this can be considered as the functional equation of an SDP, where the holding/penalty cost of period $t$ is the immediate cost. There are two boundary conditions:
\begin{align}
    l_{T+1}(I_T+1 , R_t) = 0\\
    l_{t}(I_t , 0) =  0\label{eq:boundaryl}
\end{align}

The states are represented by the tuple $(t, I_t, R_t)$ and are computed in a forward manner. To avoid recomputations, we store the computed tuples in a dictionary with constant access time.

\subsection{Unit cost and lost sales extensions} \label{sec:unit_cost}
Similarly to \cite{visentin2021computing}, unit ordering cost can be easily modelled as a function of the expected closing inventory or included in the immediate cost function. 

In the case of a stockout, the lost sales model is more common than a delay of the demand \cite{verhoef2006out}, especially in a retail setting. Lost sales models have been underrepresented in the inventory control literature \cite{bijvank2012inventory}; however, many recent works are considering mixed lost-sales and backorder configurations \cite{elhafsi2021optimal}. The model presented herein can be adapted to include partially lost sales. \cite{dos2019enhanced} defines as $\beta$ the percentage of unmet demand that is backlogged, the remaining is lost. 
The functional equation \ref{eq:recursion_rss_edit} becomes:
\begin{equation}
\label{eq:functional_eq_lost}
C^a_t(I_{t-1}) =  \underset{0 \leq Q_t \leq M \gamma_t}{\min} ( f_t(I_{t-1}, Q_t) +
E[ C^a_{t+1}(\max(I_{t-1}+ Q_t - d_t, \beta (I_{t-1}+ Q_t - d_t)) ]  ) 
\end{equation}

\section{Experimental Results}\label{sec:experimental}
This section conducts an extensive computational study of the heuristic presented in this paper. We aim to evaluate the quality of the policies computed by the heuristic and the computational effort required. In Section \ref{sec:scalability}, we assess the computational effort required to compute a policy and the quality of the policy itself under an increasing time horizon. 
An analysis of the heuristics behaviour under different demand patterns and cost parameters is presented in Section  \ref{sec:type_analysis}. Finally, we analyse an example in which the algorithm computes a near-optimal replenishment plan. 

For the experiments, we use as a comparison the branch-and-bound (BnB) technique presented in \cite{visentin2021computing}. This is the only $(R,s,S)$ solver for this problem configuration available in the literature. We use the same solver to compute the optimality gap. The solvers are:
\begin{itemize}
\item
\textbf{BnB-Guided}, the fastest branch-and-bound approach presented in \cite{visentin2021computing}. It pre-computes an initial replenishment plan using \cite{rossi2015piecewise} to improve the computational performances.
\item
\textbf{SDP}, the basic implementation of the SDP heuristic model presented in Algorithm \ref{alg:sdp}. We include this to appreciate the impact of the optimisation techniques deployed.
\item
\textbf{SDP-Opt}, the heuristic implementation deployed using the K-convexity property (Algorithm \ref{alg:sdprsskconv}) and the immediate cost memoisation.
\end{itemize}
All experiments are executed on an Intel(R) Xeon E5640 Processor (2.66GHz) with 12 Gb RAM. For the sake of reproducibility, we made the implementation of all the techniques and the data generators available\footnote{\url{https://github.com/andvise/inventory-control}}.

Since our approach is an heuristic, we use the optimality gap as measure to compute the computed policy's quality. The optimality gap is the estimated extra-cost of using the policy instead of the cost-optimal one for a particular problem. It is defined as:
\begin{equation}\label{eq:ch4optimality_gap}
    \text{Optimality gap} \triangleq \frac{\text{Policy cost} - \text{Optimal cost}}{\text{Optimal cost}}
\end{equation}
Better policy parameters exhibit a lower optimality gap. It can be used to estimate the inventory cost of deploying a non-optimal system.

\subsection{Scalability}\label{sec:scalability}
We used the same testbed presented in \cite{visentin2021computing}. A fixed holding cost per unit $h=1$. The other cost factors are sampled from uniform random variables: fixed ordering cost $K \in [80, 320]$, fixed review cost $W \in [80, 320]$ and linear penalty cost $b \in [4,16]$. The demand is modelled as a series of Poisson random variables. A uniform random variable draws the average demands per period with a range of 30 to 70. We generate 100 different instances. We replicate the experiments for increasing values of the number of periods.

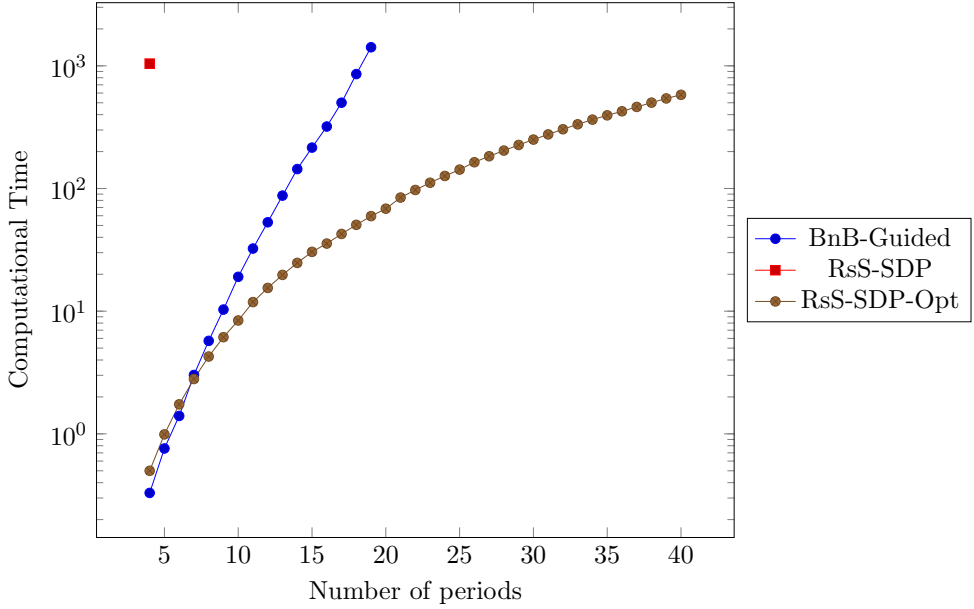
\begin{figure}
    \centering
\begin{tikzpicture}[scale = 0.9]
\pgfplotsset{every axis legend/.append style={at={(1.02,0.5)}, anchor=west}}
\begin{semilogyaxis}[xlabel={Number of periods}, ylabel={Computational Time}]

\addplot coordinates{
(4,0.33)    (5,0.76)    (6,1.4)    (7,3.02)    (8,5.73)    (9,10.31)    (10,19.06)    (11,32.41)    (12,53.05)    (13,87.54)    (14,144.08)    (15,215.53)    (16,320.05)    (17,500.68)    (18,856.64)    (19,1417.32)    
};
\addplot coordinates{
(4,1043.32)    
};

\addplot coordinates{
(4,0.5)    (5,0.99)    (6,1.74)    (7,2.8)    (8,4.27)    (9,6.14)    (10,8.4)    (11,11.86)    (12,15.47)    (13,19.76)    (14,24.75)    (15,30.45)    (16,35.58)    (17,42.62)    (18,50.58)    (19,59.54)    (20,68.5)    (21,84.39)    (22,97.27)    (23,111.42)    (24,126.58)    (25,142.59)    (26,163.71)    (27,183.24)    (28,203.81)    (29,226.43)    (30,250.5)    (31,276.08)    (32,304.11)    (33,333.73)    (34,364.4)    (35,394.84)    (36,426.0)    (37,462.46)    (38,501.19)    (39,542.28)    (40,580.47)    
};
\legend{BnB-Guided, RsS-SDP, RsS-SDP-Opt}
\end{semilogyaxis}
\end{tikzpicture}
\caption{Computational time of the $(R,s,S)$ SDP over the number of periods, time limit 1 hour}
\label{fig:scalability_log}
\end{figure}

Figure \ref{fig:scalability_log} shows the logarithm of the average computational time over the 100 instances in comparison with the fastest technique available in the literature. The simple implementation of the heuristic can barely solve tiny instances before the time limit, making it useless for every practical use. The reduction of computational effort provided by K-convexity and memoisation is massive. The guided BnB slightly outperforms the optimised SDP for small instances up to 8 periods, then the gap between the two strongly increases, making it able to solve instances more than twice as big in the same amount of time. The K-convexity performances improvement is more significant than the memoisation one. Moreover, it generally avoids the computation of all the DP states associated with a negative inventory (line 13 of Algorithm \ref{alg:sdprsskconv}). The memoisation offers a great speed up in the computational times, which is more significant in bigger instances. For bigger instances, the physical memory needed grows to require the usage of memory swap and a slow down in performances.

In this testbed, the heuristic always computes the optimal replenishment plan.

\subsection{Instance type analysis}
\label{sec:type_analysis}
These experiments aim to analyse the performances of the heuristic under different instance parameters. We want to analyse which cost parameters are affecting the computational performances and the optimality gap of the heuristic.
We use a modified version of the instances used in Section 6.2 of \cite{visentin2021computing}. The algorithm proposed herein computes the optimal policy parameters for all the instances used therein. Our extension aims to find problem settings where the heuristic under-performs the optimal approach. We do it by examining a wider range of review and ordering costs and increasing the demand's uncertainty. 
Poisson distributed demand does not have a parameter to increase the uncertainty over its expected value; for this reason, we included normal demand in our experiments.

We use two different planning horizons: $10$ and $20$ periods. For the cost parameters, we use all the possible combinations of review and ordering cost values $K, W \in \{20 , 40, 80,160,320\}$, holding and penalty cost fixed respectively at $h = 1$ and $b = 10$. We consider Poisson demand and normally distributed one with  $\sigma \in \{0.1, 0.2, 0.3, 0.4\}$. In the literature, the standard deviation used is generally not higher than $0.3$; we use $0.4$ to simulate extreme uncertain cases. We consider six different demand patterns: stationary (STA), positive trend (INC), negative trend (DEC), two life-cycle trends (LCY1, LCY2) and an erratic one (RAND); more details on these patterns can be found in \cite{visentin2021computing}. The combinations of the parameters mentioned above lead to the creation of 1 500 instances.

Table  \ref{tab:10_periods}  and Table \ref{tab:20_periods} show the results for the 10 and 20-period instances. Regarding the policies' quality, we consider the average optimality gap, the percentage of computed policies that differs for the optimal and their optimality gap. We also compare the time required to compute the policies and the average number of reviews.

The cost factors suggest that the algorithm does not compute the optimal policy in situations with a high ordering cost and a low review cost. Due to the relaxation, the approximate SDP computes the parameters of the cycles based only on the state values, considering the uncertainty of the future periods but ignoring the one related to the period opening inventory. When the review cost is low, the BnB uses more review periods compared to the SDP to counteract this uncertainty. Up to $16.67\%$ and  $18.67\%$ (for the 10 and 20 instances) of the policy computed differ from the optimal one; however, their gap averages less than $1\%$. The average optimality gap across all the instances with the lowest review cost is $0.22\%$ and $0.29\%$. A higher ordering cost leads to longer intervals between orders, so a higher uncertainty on the opening inventory level of a period. This leads to a maximum of $10.67\%$ and $13.33\%$ of near-optimal policies. While for these particular settings, the percentage of non-optimal policy is relatively high; their optimality gap is low.

The direct correlation between the demand uncertainty and the optimality gap is evident. In the literature, the standard deviation used is in the range $[0.1, 0.3]$; we used expected demand with a higher degree of uncertainty to show more clearly the situations in which our approach struggles. With realisations of the demand that strongly differ from their expected value, it is more likely that the opening inventory level is higher than the order-up-to-level in a review moment. In these cases, our approach relaxes the problem by placing a negative order and setting the inventory to $S$, so the policy is not optimal. The majority of the instances in which the SDP computes a near-optimal policy have a $\sigma = 0.4$.

The pattern analysis provides interesting insights. The approach performs better with the increasing (INC) pattern regardless of the other instance parameters; in 10-periods instances, it always computes the optimal policy. We have the worse performances in the decreasing (DEC) one. This is in line with \cite{ozen2012static} that considers a similar problem relaxation. If we have increasing demands, the base stock levels likely increase as well to satisfy higher demands. If the base stock levels increase monotonically, the relaxation generally computes the optimal policy. The second worse pattern is the random one due to randomly generated decreasing patterns. We observe the biggest gap between the number of reviews with $0.2$ and $0.5$ fewer reviews on average in the decreasing pattern.

On average, the optimality gap between the two approaches is only $0.04\%$ and  $0.05\%$, with $4.4\%$ and $5.87\%$ of the policy computed that are near-optimal. This proves the quality of the heuristic in computing policies.

Our approach is 4 and 300 times faster, respectively, on the 10 and 20 periods regarding the computational time. Moreover, the cost parameters do not affect the SDP performances, while they affect the BnB pruning efficacy. For example, low review cost 20-period instances takes six times more effort than high review ones. Uncertainty on the forecast affects the performances of both approaches since it makes the computation of a state expected cost more expensive. However, the SDP manages to reduce this impact using memoisation. The SDP increases less than  times its computational effort for $\sigma = 0.4$ compared to $\sigma = 0.1$, while the increment for the BnB approach is higher than 40.

\begin{table*}
\centering
\resizebox{\hsize}{!}{%
\def\arraystretch{1.4}
\renewcommand{\tabcolsep}{1.5mm}
\begin{tabular}{|c c||c|c|c||c|c||c|c||c|c|}
\hline
& & \multicolumn{3}{c||}{Optimality} & \multicolumn{2}{c|}{Time (min)} & \multicolumn{2}{c|}{Nr Reviews} & \multicolumn{2}{c|}{Expected cost error}  \\ \hline
& & Optimality Gap &  \% Non-Optimal  & Non-Optimal OG & SDP & BnB & SDP & BnB & SDP & BnB    \\ \hline
K values & 20 & 0.0	& 0.0	& 0.0	& 0.27	& 0.94	& 5.28	& 5.29	& 0.2	& 0.2\\ 
		 & 40 & 0.0	& 1.33	& 0.26	& 0.27	& 1.04	& 4.63	& 4.65	& 0.2	& 0.19\\ 
		 & 80 & 0.01	& 4.67	& 0.35	& 0.28	& 1.13	& 3.77	& 3.8	& 0.17	& 0.17\\ 
		 & 160 & 0.03	& 5.33	& 0.67	& 0.28	& 1.3	& 3.13	& 3.21	& 0.15	& 0.14\\ 
		 & 320 & 0.1	& 10.67	& 1.07	& 0.3	& 1.49	& 2.42	& 2.56	& 0.11	& 0.11\\ 
\hline
W values & 20 & 0.22	& 16.67	& 0.97	& 0.26	& 1.41	& 5.37	& 5.59	& 0.19	& 0.18\\ 
		 & 40 & 0.03	& 4.0	& 0.58	& 0.27	& 1.44	& 4.64	& 4.69	& 0.18	& 0.18\\ 
		 & 80 & 0.0	& 0.67	& 0.33	& 0.28	& 1.34	& 3.77	& 3.77	& 0.17	& 0.17\\ 
		 & 160 & 0.0	& 0.67	& 0.07	& 0.29	& 1.02	& 3.09	& 3.09	& 0.15	& 0.15\\ 
		 & 320 & 0.0	& 0.0	& 0.0	& 0.3	& 0.68	& 2.37	& 2.37	& 0.13	& 0.13\\ 
\hline
Poisson & p & 0.0	& 0.0	& 0.0	& 0.17	& 0.47	& 3.65	& 3.65	& 0.12	& 0.12\\ 
$\sigma$ values & 0.1 & 0.0	& 0.67	& 0.07	& 0.16	& 0.24	& 3.68	& 3.68	& 0.08	& 0.08\\ 
		 & 0.2 & 0.01	& 1.33	& 0.66	& 0.24	& 0.69	& 3.81	& 3.81	& 0.16	& 0.16\\ 
		 & 0.3 & 0.05	& 7.33	& 0.66	& 0.36	& 1.74	& 3.95	& 4.05	& 0.2	& 0.2\\ 
		 & 0.4 & 0.12	& 12.67	& 1.01	& 0.47	& 2.76	& 4.14	& 4.31	& 0.26	& 0.24\\ 
\hline
Pattern & STA & 0.03	& 4.8	& 0.73	& 0.15	& 0.97	& 3.87	& 3.95	& 0.15	& 0.13\\ 
		 & INC & 0.0	& 0.0	& 0.0	& 0.37	& 1.17	& 4.05	& 4.05	& 0.15	& 0.15\\ 
		 & DEC & 0.11	& 9.6	& 1.12	& 0.33	& 1.73	& 3.47	& 3.56	& 0.14	& 0.14\\ 
		 & LCY1 & 0.02	& 1.6	& 1.2	& 0.24	& 1.09	& 4.02	& 4.07	& 0.2	& 0.2\\ 
		 & LCY2 & 0.02	& 4.0	& 0.44	& 0.3	& 1.14	& 3.86	& 3.91	& 0.2	& 0.2\\ 
		 & ERR & 0.05	& 6.4	& 0.71	& 0.29	& 0.97	& 3.8	& 3.87	& 0.14	& 0.14\\ 
\hline
Average & & 0.04	& 4.4	& 0.84	& 0.28	& 1.18	& 3.85	& 3.9	& 0.16	& 0.16\\ 
\hline
\end{tabular}
}
\caption{Optimality gap and pruning percentage for the techniques for instances of 10 periods}
\label{tab:10_periods}
\end{table*}

\begin{table*}
\centering
\resizebox{\hsize}{!}{%
\def\arraystretch{1.4}
\renewcommand{\tabcolsep}{1.5mm}
\begin{tabular}{|c c||c|c|c||c|c||c|c||c|c|}
\hline
& & \multicolumn{3}{c||}{Optimality} & \multicolumn{2}{c|}{Time (min)} & \multicolumn{2}{c|}{Nr Reviews} & \multicolumn{2}{c|}{Expected cost error}  \\ \hline
& & Optimality Gap &  \% Non-Optimal  & Non-Optimal OG & SDP & BnB & SDP & BnB & SDP & BnB   \\ \hline
K values & 20 & 0.0	& 0.0	& 0.0	& 2.2	& 398.65	& 10.47	& 10.49	& 0.05	& 0.05\\ 
		 & 40 & 0.01	& 3.33	& 0.21	& 2.12	& 487.29	& 9.21	& 9.29	& 0.05	& 0.05\\ 
		 & 80 & 0.01	& 5.33	& 0.21	& 2.11	& 639.83	& 7.45	& 7.51	& 0.05	& 0.05\\ 
		 & 160 & 0.05	& 7.33	& 0.72	& 2.23	& 848.22	& 6.15	& 6.33	& 0.05	& 0.05\\ 
		 & 320 & 0.13	& 13.33	& 1.06	& 2.05	& 1033.47	& 4.71	& 5.06	& 0.05	& 0.05\\ 
\hline
W values & 20 & 0.29	& 21.33	& 0.99	& 2.01	& 1050.88	& 10.65	& 11.17	& 0.05	& 0.05\\ 
		 & 40 & 0.05	& 6.67	& 0.45	& 2.14	& 988.62	& 9.24	& 9.4	& 0.05	& 0.05\\ 
		 & 80 & 0.0	& 1.33	& 0.1	& 2.13	& 787.85	& 7.41	& 7.41	& 0.05	& 0.05\\ 
		 & 160 & 0.0	& 0.0	& 0.0	& 2.18	& 418.7	& 6.07	& 6.07	& 0.05	& 0.05\\ 
		 & 320 & 0.0	& 0.0	& 0.0	& 2.24	& 161.41	& 4.62	& 4.62	& 0.05	& 0.05\\ 
\hline
Poisson & p & 0.0	& 0.0	& 0.0	& 1.34	& 130.59	& 7.37	& 7.37	& 0.04	& 0.04\\ 
$\sigma$ values & 0.1 & 0.0	& 0.0	& 0.0	& 1.27	& 46.83	& 7.37	& 7.37	& 0.03	& 0.03\\ 
		 & 0.2 & 0.0	& 1.33	& 0.14	& 1.89	& 279.82	& 7.51	& 7.51	& 0.05	& 0.05\\ 
		 & 0.3 & 0.05	& 9.33	& 0.56	& 2.69	& 946.02	& 7.71	& 7.9	& 0.06	& 0.06\\ 
		 & 0.4 & 0.17	& 18.67	& 0.98	& 3.52	& 2004.2	& 8.03	& 8.53	& 0.07	& 0.07\\ 
\hline
Pattern & STA & 0.05	& 4.8	& 1.07	& 0.95	& 574.5	& 7.57	& 7.8	& 0.02	& 0.02\\ 
		 & INC & 0.01	& 2.4	& 0.51	& 2.6	& 347.01	& 7.74	& 7.77	& 0.05	& 0.05\\ 
		 & DEC & 0.08	& 5.6	& 1.36	& 2.52	& 1438.91	& 7.15	& 7.33	& 0.09	& 0.08\\ 
		 & LCY1 & 0.04	& 4.8	& 0.7	& 1.74	& 481.35	& 7.98	& 8.09	& 0.05	& 0.05\\ 
		 & LCY2 & 0.04	& 8.0	& 0.49	& 2.4	& 517.99	& 7.48	& 7.6	& 0.05	& 0.05\\ 
		 & ERR & 0.08	& 9.6	& 0.78	& 2.63	& 729.19	& 7.67	& 7.83	& 0.05	& 0.05\\ 
\hline
Average & & 0.05	& 5.87	& 0.81	& 2.14	& 681.49	& 7.6	& 7.74	& 0.05	& 0.05\\ 
\hline
\end{tabular}
}
\caption{Optimality gap and pruning percentage for the techniques for instances of 20 periods}
\label{tab:20_periods}
\end{table*}

\subsubsection{Non-optimality of the relaxation}
In this section, we analyse a single instance to better understand the differences between the computed policies. This example shows a situation in which the heuristic computes a non-optimal policy. When computing the solution, it considers only the expected demand for future periods. On the other hand, the BnB approach presented in \cite{visentin2021computing} tests all the possible replenishment combinations of the previous periods during the search process. Not considering the previous demands means ignoring the possibility of having such a low demand that at a period $t$, the opening inventory level $I_t$ is higher than $s_t$, and that this will strongly affect future decisions. This difference worsens the heuristics performances for high values of uncertainty and the decreasing pattern (DEC). In these instances, the high demand with high uncertainty at the beginning of the time horizon makes unexpected high inventory levels at a replenishment moment more likely. In this situation, the BnB solution adds more review moments (especially when the cost associated $W$ is low) to assess the inventory level and react to the uncertainty.

For example, considering the instance of Table \ref{tab:10_periods} with $K = 320$, $W=20$, $\sigma = 0.4$ and decreasing demand pattern. Table \ref{tab:policycomparison} shows the two policies. The BnB approach considers the higher uncertainty at the beginning of the time horizon; it also reviews the inventory level at periods $5$ and $6$. While these reviews add an extra cost in an almost deterministic system, they allow a better reaction to unexpected demand. At the end of the time horizon, the uncertainty on the inventory level is lower, and the two policies are identical from period $7$ on when a lower demand leads to lower absolute variations of the realised demand.
 
The BnB policy has an expected cost of 1793, the SDP of 1845; a difference of 52 that leads to an optimality gap of $2.9\%$.

\begin{table*}
\centering
\resizebox{0.9\textwidth}{!}{%
\def\arraystretch{1.5}
\renewcommand{\tabcolsep}{1.5mm}
\begin{tabular}{|c x{0.8cm}||x{0.8cm}|x{0.8cm}|x{0.8cm}|x{0.8cm}|x{0.8cm}|x{0.8cm}|x{0.8cm}|x{0.8cm}|x{0.8cm}|x{0.8cm}||c|}
\hline
 \multicolumn{2}{|c||}{Period}&1 &2 &3 &4 &5 &6 &7 &8 &9 &10& Policy cost\\
\hline
\hline
        &$\gamma_t$ & 1& 0& 0& 1& 0& 0& 0& 1&0& 0&\\
RsS-SDP &$S_t$ &295& -& -& 243& -& -& -& 56& -& -&1845\\
        &$s_t$ &211& -& -& 174& -& -& -& 25& -& -&\\
\hline
\hline
        &$\gamma_t$ &1& 0& 0& 1& 1& 1& 0& 1& 0& 0&\\
RsS-BnB &$S_t$ &324& -& -& 237& 186& 139& -& 56& -& -&1793\\
        &$s_t$ &220& -& -& 48& 42& 64& -& 25& -& -&\\
        \hline
\end{tabular}
}
\caption{$(R,s,S)$ policy parameters for the  $K = 320$, $W=20$, $\sigma = 0.4$, DEC pattern instance.}
\label{tab:policycomparison}
\end{table*}

\section{Conclusions} \label{sec:conclusions}
This paper presented a heuristic for the non-stationary stochastic lot-sizing problem with ordering, review, holding and penalty cost, a well-known and widely used inventory control problem. Computing $(R,s,S)$ policy parameters is computationally hard due to the three sets of parameters that must be jointly optimised. We presented the first pure SDP formulation for such a problem. The algorithm introduced solves to optimality a relaxation of the original problem, in which review cycles are considered independently, and items can be returned/discarded at no additional cost. A similar relaxation has been previously used in $(R,S)$ policy computation works.

The extensive numerical study proved the reduction of the computational effort needed to compute a policy. The basic formulation requires a prohibitive computational effort. Two enhancements based on K-convexity \cite{scarf1959optimality} and memoisation strongly improve the computational performance, making it able to solve instances twice as big as the state-of-the-art. This allows practitioners to use such policy in a wider range of real-world situations. We then investigated the SDP performance under different types of instances. We measured the computational effort to compute the policy and how much the relaxation affects its quality. The heuristics' computational effort is less affected by the instance configuration. The proposed algorithm rarely computes a non-optimal policy when there is less uncertainty on demand and high review, low fixed ordering cost instances. For Poisson distributed demand, the SDP always computes the optimal policy. The average optimality gap is $0.04 \%$ and $0.05\%$ with $95.6\%$ and $94.13\%$ of computed policies identical to the optimal respectively for the 10 and the 20 periods instances; more than half of the non-optimal policies are related to extremely high uncertainty of the demand ($\sigma = 0.4$) a configuration hardly considered in the lot-sizing literature. These differences are caused by a reduced number of review moments in the SDP computed policies.  

In future studies, we plan to extend such a method's applicability by considering more complex supply chains such as multiple items, multiple echelons, and different cost structures. We plan to further enhance the current formulation to improve the non-optimal computed policies, similarly to what \cite{rossi2011state} did with state space augmentation.

\subsubsection*{Acknowledgments}
This publication has emanated from research conducted with the financial support of Science Foundation Ireland under Grant number 16/RC/3918  which is co-funded under the European Regional Development Fund. For the purpose of Open Access, the author has applied a CC BY public copyright licence to any Author Accepted Manuscript version arising from this submission.

\bibliographystyle{tfcad}
\bibliography{mybibfile}

\end{document}